# The Derivative of the Sine and Cosine. A New Derivation Approach

*By John T. Katsikadelis*

Scool of Civil Engineering, National Technical University of Athens, Athens 15773, Greece.

e-mail: jkats@central.ntua.gr    Website: http://users.ntua.gr/jkats/

**Abstract**

A new method is presented for finding the derivative of the sine and cosine using the discoveries of Leibniz in calculus between the years 1675 and 1677, namely the derivative of the product and the quotient of two functions as well as the chain rule, yet long before the discovery of the derivative of the sine by Roger Cotes in 1722 and Euler in 1748.

**Introduction**

The seeds for finding of the derivative of the sine and cosine have been sown by Roger Cotes (1682-1716) in his work "*Harmonia mensurarum*" [1] published posthumously in 1722. Cotes stated and proved the following lemma:

> **Lemma I.**
>
> **Variatio minima cujusvis arcus circularis est ad Variationem minimam Sinus ejusdem arcus ut radius ad Sinum complementi.**

which is rendered in English as [2]

> **Lemma I.**
>
> **The small variation of any arc of a circle is to the small variation of the sine of that arc, as the radius to the sine of the complement.**

The proof given by Cotes results from Fig. 1 (quoted from [1]) by considering the similar triangles ACD and ECG. Thus we have

$$\frac{EC}{EG} = \frac{AC}{AD} \tag{1}$$

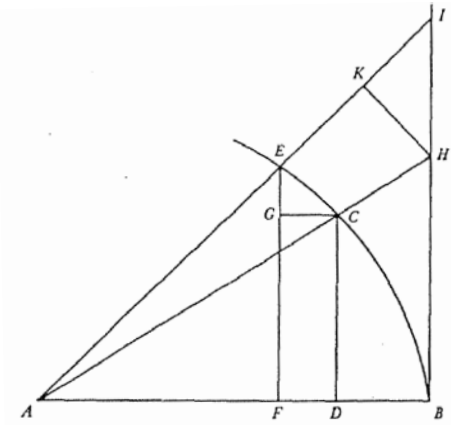

Figure 1.

Using modern notation, the *small variation of any arc of a circle* with radius $r$ becomes $d(r\theta) = CE$, while the *small variation of the sine of that arc* is $d(\sin\theta) = CG$. Hence, Eq. (1) is written as

$$\frac{d(r\theta)}{d(\sin\theta)} = \frac{r}{\cos\theta} \qquad (2)$$

In 1748, Euler introduced the trigonometric functions to treat the linear homogeneous differential equations [3] and simplified Eq. (2) using the unit circle, yielding thus the derivative of the sine

$$\frac{d(\sin\theta)}{d\theta} = \cos\theta \qquad (3)$$

Since the proof given by Cotes might be objected as lacking in rigor, the derivative was obtained later using the relation derived by Cauchy in 1822 [4]

$$\lim_{\theta \to 0} \frac{\sin\theta}{\theta} = 1 \qquad (4)$$

in which the measure of the angles is in radians. Eq. (4) results for small values of $\theta$ from the inequality

$$1 > \frac{\sin\theta}{\theta} > \cos\theta \qquad (5)$$

which had been implied in the works of Archimedes. Hence we obtain

$$\frac{d}{d\theta}(\sin\theta) = \lim_{\Delta\theta \to 0} \frac{\sin(\theta + \Delta\theta) - \sin\theta}{\Delta\theta}$$

$$= \lim_{\Delta\theta \to 0} \left[\cos(\theta + \Delta\theta / 2) \frac{\sin(\Delta\theta / 2)}{\Delta\theta / 2}\right] \quad (6)$$

$$= \lim_{\Delta\theta \to 0} \cos(\theta + \Delta\theta / 2) \lim_{\Delta\theta \to 0} \frac{\sin(\Delta\theta / 2)}{\Delta\theta / 2}$$

$$= \cos\theta$$

This author in his recent paper [5] avoided the differentiation of the cosine because he wanted to use mathematics that was available before 1686, namely the year when Newton published his Principia [6]. But this was still a challenge to develop a new method for finding the derivative of the sine and cosine using mathematics available before 1686. The presented procedure uses the rules for the differentiation of the product and the quotient of two functions of one variable as well as the chain rule. All these rules were discovered by Leibniz about 1676 [7].

**The new approach for finding the derivative of the sine and cosine**

We consider a moving particle in the $xy$ plane. Its position in polar coordinates determined by the time dependent parameters $r(t)$ and $\theta(t)$, while in Cartesian coordinates by $x(t)$ and $y(t)$, Fig. 2.

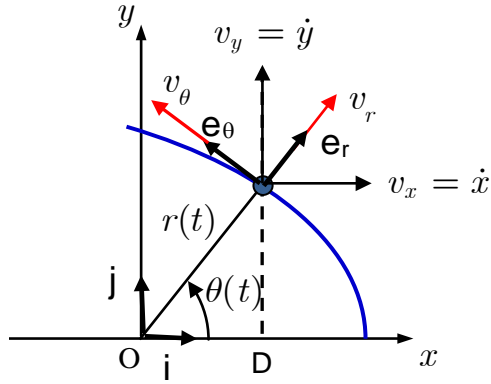

Figure 2

We define the functions:

$$c = \cos\theta = \frac{x}{r}, \quad s = \sin\theta = \frac{y}{r} \quad (7)$$

Then we have

$$x = rc, \qquad y = rs, \qquad r^2 = x^2 + y^2 \tag{8a,b,c}$$

Differentiation of Eqs. (8) with respect to $t$ gives

$$\dot{x} = \dot{r}c + r\dot{c} \tag{9a}$$

$$\dot{y} = \dot{r}s + r\dot{s} \tag{9b}$$

$$\dot{r} = \frac{x\dot{x} + y\dot{y}}{r} \tag{9c}$$

Multiplying Eqs. (9a) and (9b) by $y$ and $x$, respectively, and subtracting them, we obtain

$$y\dot{c} - x\dot{s} = \frac{y\dot{x} - x\dot{y}}{r} \tag{10}$$

Further, multiplying Eqs. (9a) and (9b) by $x$ and $y$, respectively, adding them and using Eq. (8c), we obtain

$$x\dot{x} + y\dot{y} = r\dot{r} + rx\dot{c} + ry\dot{s} \tag{11}$$

which by virtue of Eq. (9c) becomes

$$x\dot{c} + y\dot{s} = 0 \tag{12}$$

Eqs. (10) and (12) are combined and solved for $\dot{s}$ and $\dot{c}$. Thus we obtain

$$\dot{s} = -\frac{x}{r}\frac{y\dot{x} - x\dot{y}}{r}\frac{1}{r} = -\cos\theta(\dot{x}\sin\theta - \dot{y}\cos\theta)\frac{1}{r} \tag{13a}$$

$$\dot{c} = \frac{y}{r}\frac{y\dot{x} - x\dot{y}}{r}\frac{1}{r} = \sin\theta(\dot{x}\sin\theta - \dot{y}\cos\theta)\frac{1}{r} \tag{13b}$$

Referring to the Fig. 2, we have

$$\dot{y}\cos\theta - \dot{x}\sin\theta = v_\theta = r\dot{\theta} \tag{14}$$

where $v_\theta$ is the component of the velocity in the direction $\mathbf{e}_\theta$, i.e., normal to $\mathbf{e}_r$.

Hence, Eqs. (13a,b) by virtue of Eq. (14) and the use of the chain rule become

$$\dot{s} = \frac{d\sin\theta}{d\theta}\dot{\theta} = \cos\theta\dot{\theta} \tag{16a}$$

$$\dot{c} = \frac{d\cos\theta}{d\theta}\dot{\theta} = -\sin\theta\dot{\theta} \tag{16b}$$

or

$$\frac{d\sin\theta}{d\theta} = \cos\theta \tag{17a}$$

$$\frac{d\cos\theta}{d\theta} = -\sin\theta \tag{17b}$$

Obviously, Eqs. (17a,b) give the derivative of the sine and cosine with respect to the angle measured in radians.

**Conclusions**

In this note a new method is presented for finding the derivative of the sine and cosine using mathematics available before 1686, the year Newton published his Principia, namely the derivative of the product and the quotient of two functions as well as the chain rule discovered by Leibniz between the years 1675 and 1677.